\documentclass[a4paper]{article}
\usepackage[english]{babel}
\usepackage[T1]{fontenc}
\usepackage{lmodern}
\usepackage{chapterbib}
\usepackage{amsthm}
\usepackage{siunitx}
\usepackage{subcaption}
\usepackage{graphicx}
\usepackage{xspace}
\usepackage[inkscapeformat=pdf]{svg}
\usepackage{amsmath,amsfonts,amssymb,amsthm}

\usepackage[colorlinks=true,allcolors=blue]{hyperref}
\usepackage[a4paper,top=3cm,bottom=2.2cm,left=2.65cm,right=2.65cm]{geometry}

\theoremstyle{remark}

\theoremstyle{definition}

\usepackage{graphicx}
\usepackage{amsmath}
\usepackage{natbib}
\usepackage{float}
\usepackage{authblk}
\usepackage{booktabs}

\usepackage{subcaption}
\usepackage{multicol}
\usepackage{array}
\usepackage{arydshln}
\newcommand\Tstrut{\rule{0pt}{2.6ex}}         
\newcommand\Bstrut{\rule[-0.9ex]{0pt}{0pt}}   

\usepackage{xcolor}
\usepackage{tikz}
\newcommand{\x}{\boldsymbol{x}}
\providecommand{\abs}[1]{\lvert#1\rvert}

\newcommand{\xb}{\boldsymbol{\bar{x}}}

\newcommand{\Xb}
{\boldsymbol{\bar{X}}}
\newcommand{\R}{\mathbb{R}}

\newcommand{\N}{\mathbb{N}}

\newcommand{\nf}[1]{\footnotesize{#1}}

\usepackage{todonotes}

\title{An extension of an RLT-based solver to MINLP polynomial problems} 
\author[1,2]{Julio González-Díaz}
\author[1]{Brais González-Rodríguez}
\author[1]{Iria Rodríguez-Acevedo}
\affil[1]{Department of Statistics, Mathematical Analysis and Optimization and MODESTYA Research Group,  University of Santiago de Compostela.}
\affil[2]{CITMAga (Galician Center for Mathematical Research and Technology).}

\date{\today}

\newcommand{\solver}[1]{\texttt{#1}}
\newcommand{\CC}{C\nolinebreak\hspace{-.05em}\raisebox{.4ex}{\tiny\bf +}\nolinebreak\hspace{-.10em}\raisebox{.4ex}{\tiny\bf +}\xspace}

\usepackage{mathtools}

\usepackage[skins]{tcolorbox}

\usepackage{enumitem}
\setlist[itemize,3]{label=$\diamond$}

\newenvironment{shortitem}[1][]{\begin{itemize}[topsep=0pt, itemsep=0pt, parsep=0.5pt, leftmargin=*, #1]}{\end{itemize}}

\newcommand{\itercounter}{\tau}

\newcommand{\cola}{Q}

\newcommand{\zbar}{\bar z}

\renewcommand{\stop}{\textsc{End}}

\begin{document}

\maketitle

\begin{abstract}
In this paper we extend the core branch-and-bound algorithm of an RLT-based solver for continuous polynomial optimization, \solver{RAPOSa}, to handle mixed-integer problems. We do so by a direct adaptation, in which LP relaxations are replaced with MILP ones and, therefore, the additional burden caused by the discrete variables is taken care of by the auxiliary MILP solver. The analysis then focuses on the impact of a number of choices that must be made for the resulting algorithm to be effective at providing good lower and upper bounds upon termination. 
\end{abstract}

\section{Introduction}

In the recent years there has been an explosion in the design and implementation of global optimization algorithms for MINLP problems. State-of-the-art commercial solvers for MILP problems such as \solver{Gurobi} \citep{gurobi} and \solver{Xpress} \citep{xpress} have recently expanded their functionalities to solve MINLP problems to certified global optimality. \solver{Octeract} \citep{octeract} is another global solver that has recently joined the list of well-established solvers, which includes \solver{ANTIGONE} \citep{antigone}, \solver{BARON} \citep{Khajavirad2018,baron}, \solver{Couenne} \citep{Belotti2009}, \solver{LindoGlobal} \citep{lindo}, and \solver{SCIP} \citep{scip}. This rise of global optimization algorithms is being fueled by the increase in the processing power of the computers on which they can be run and by the generality of MINLP problems, which leads to a plethora of applications in a wide variety of fields.

All of the above solvers rely on different transformations and reformulations of the original MINLP problems, with the goal of obtaining easy-to-solve relaxations that can be embedded into branch-and-bound schemes.\footnote{We refer the reader to \cite{Burer2012} and \cite{Belotti2013} for two complete surveys on the elements involved in the design of global optimization algorithms for MINLP problems.} Undoubtedly, one of the most important reformulations is the so-called reformulation-linearization technique, hereafter RLT, central to this paper's main contribution: to enhance the RLT-based solver \solver{RAPOSa} \citep{GonzalezRodriguez2023}, expanding its capabilities from continuous polynomial optimization problems to also address mixed-integer ones.

RLT was formally introduced in \cite{SheraliAdams1990,SheraliAdams1994}, where the authors present a hierarchy of relaxations of pure 0-1 and mixed-integer 0-1 linear and polynomial programming problems.\footnote{Actually, two slightly older papers, \cite{SheraliAdams1984}  and \cite{AdamsSherali1986}, already included RLT-type liftings to generate tighter representations of some specific problems, but the name RLT was coined in \cite{SheraliAdams1990,SheraliAdams1994}. Refer to \cite{Sherali2007} for a historical review of the origins of the technique.} A few years later, \cite{SheraliAdams1998} strengthened the results in \cite{SheraliAdams1994} and \citeauthor{SheraliAdams1999} (\citeyear[Chapter~4]{SheraliAdams1999}) extended the analysis to general mixed-integer linear problems. The main goal of these series of papers was to provide hierarchies of sharper and sharper relaxations, starting with the continuous one and ending with a final relaxation representing the convex hull of the set of feasible solutions. Thus, solving this final relaxation would be equivalent to solving the original MINLP problem. 

Since the number of variables in the RLT relaxations increases as one moves up in the hierarchy, one of the most important challenges of RLT is to keep the size of the problems within manageable limits. Alternatively, there is significant research on approaches that start with looser relaxations and then focus on finding classes of valid inequalities which, eventually, might also end up with the convex hull of the set of feasible solutions \citep{crama2017,delpiaKhajavirad2018,delpiaKhajavirad2021,delpiaSahinidis2020,delpiaGregorio2021,delpiaGregorio2023,delpiaWalter2023}. An intermediate (and complementary) approach can be seen in \cite{elloumi2023}, where the authors present a formal model for ``extended linear reformulations of polynomials of binary variables'' with the goal of optimizing the choice of these linear patterns, either to maximize the dual bound or to minimize the size of the resulting relaxations.

Importantly, in a relatively independent series of papers, RLT was also applied to the solution of continuous polynomial optimization problems. \cite{Sherali1992} present a global optimization algorithm to solve these problems based on RLT relaxations, and subsequent studies aimed at improving these relaxations, either by reducing their size \citep{Sherali1997,Sherali2012reduced,Dalkiran2013} or improving their tightness \citep{Sherali2011}. \solver{RLT-POS} \citep{Dalkiran2016} and \solver{RAPOSa} \citep{GonzalezRodriguez2023} are two RLT-based solvers for continuous polynomial optimization problems.

Naturally, there are tight connections between the use of RLT for discrete and continuous optimization problems and, yet, we are not aware of any research or RLT-based solver integrating them. This observation is what drives the main contribution of this paper, which takes a first step in this direction by presenting a first implementation of an RLT-based solver that can address mixed-integer polynomial optimization problems, while certifying global optimality. We do so by means of a direct adaptation of the core algorithm in the global solver \solver{RAPOSa}, which has proven very effective in our preliminary analysis. Namely, the LP relaxations used to solve continuous problems are now replaced by MILP ones and, therefore, the integrality is handled by the auxiliary MILP solver. Thus, full priority is given to branching on integer variables, which agrees with the approach taken in the default configuration of \solver{Couenne} \citep{Belotti2009}. \solver{Xpress} \citep{xpress}, on the other hand, determines the branching variable by computing scores along the lines of reliability branching \citep{Achterberg2005}, with no particular priority given to either continuous or integer variables. \cite{kilinc2018} and \cite{zhou2018} present a dynamic strategy for deciding under what conditions to activate integer programming relaxations in the course of the branch-and-bound scheme in \solver{BARON} \citep{baron}, which, to some extent, would be as dynamically deciding to give priority to integer variables at certain nodes.

Despite of the simplicity of the above approach to obtain an RLT-based MINLP solver, there are a number of decisions that have to be made and a number of trade-offs that have to be properly balanced. These decisions include the need to complement the calls to an auxiliary NLP solver with calls to a computationally more expensive MINLP local solver. The speed of convergence of the upper bound of the branch-and-bound algorithm may crucially depend on adequately coordinating these two auxiliary solvers. Other decisions involve the adaptation of the domain reduction techniques for the discrete variables or reassessing the default branching rule used by the algorithm.

To the best of our knowledge, no past research has discussed the potential impact of basic decisions like the ones above in the performance of the resulting algorithm. Thus, the main contribution of the numerical analyses in this paper is to shed some light on how to properly balance the different trade-offs of such decisions. Since they are not specific to polynomial optimization, we hope that our insights can guide the design and implementation of other global optimization solvers for general MINLP problems. One of our main findings is that, in the particular case of the RLT-based algorithm under study, by carefully making the above decisions one can dramatically reduce the number of problems in which it terminates without finding either a lower or an upper bound. Importantly, we also show that, when these findings are incorporated into the freely available solver \solver{RAPOSa}, its performance becomes highly competitive with that of two state-of-the-art solvers: \solver{BARON} \citep{baron} and \solver{Couenne} \citep{Belotti2009}. Given these promising results, obtained with a relatively straightforward adaptation of the original continuous RLT-based algorithm, we believe that there is a significant potential of improvement by building upon the extensive literature on reformulations of binary and discrete variables described above, which should definitely be the main avenue for future research.

The rest of the paper is structured as follows. In Section~\ref{sec:prelim} we present some RLT preliminaries and also the setup for the computational experiments. In Section~\ref{sec:RLTint} we present a first extension of the continuous RLT-based algorithm to the mixed-integer setting. Section~\ref{sec:bounds} is devoted to study different approaches to improve the quality of the bounds obtained by the algorithm and to present a comparison of the final algorithm implemented in \solver{RAPOSa} with other state-of-the-art solvers. We conclude in Section~\ref{sec:conclusions}.

\section{Preliminaries: RLT and computational setup}\label{sec:prelim}

Given that the main object of study in this paper is the extension of a continuous RLT-based algorithm to mixed-integer problems, we start by providing a concise overview of RLT, as introduced in \cite{Sherali1992}. A minimum understanding of this technique is needed for the comprehension of the ensuing analysis, whose computational setup is described at the end of this section.

\subsection{RLT for continuous polynomial optimization }\label{sec:RLT}
In \cite{Sherali1992}, RLT was applied to solve polynomial optimization problems of the following form:
\begin{equation}
\begin{split}
\text{minimize} & \quad \phi_0(\mathbf{x})\\
\text{subject to}  & \quad \phi_r(\mathbf{x})\geq \beta_r, \quad r=1,\ldots, R_1 \\
& \quad \phi_r(\mathbf{x})=\beta_r, \quad r=R_1+1,\ldots,R\\
& \quad \mathbf{x}\in\Omega \subset \mathbb{R}^{|N|}\text{,}
\end{split}
\tag*{$PP(\Omega)$}
\label{eq:PO}
\end{equation}
where $N$ denotes the set of variables, each $\phi_r(\mathbf{x})$ is a polynomial of degree $\delta_r \in \mathbb{N}$, and $\Omega = \lbrace \mathbf{x} \in \mathbb{R}^{|N|}: 0 \leq l_j \leq x_j \leq u_j < \infty, \, \forall j \in N \rbrace \subset \mathbb{R}^{|N|}$ is a hyperrectangle containing the feasible region.\footnote{Nonnegativity of the $x_j$ variables can assumed without loss of generality, and facilitates the exposition.} Then, $\delta=\max_{r \in \{0,\ldots,R\}} \delta_r$ is the degree of the problem and $N^\delta$ represents all possible monomials of degree $\delta$.

RLT builds upon successive linearizations of problems of the form \ref{eq:PO}. The linearization process entails the following steps. First, given a monomial $J\subset N^\delta$ with $2\leq\abs{J}\leq\delta$, we define the 
corresponding RLT identity as 
\begin{equation}
  X_J = \prod_{j \in J}x_j,
  \label{eq:RLTidentity}
\end{equation}
where each $X_J$ denotes an RLT variable. Then, we define the linearization $[\phi(\x)]_L$ of a given polynomial $\phi(\x)$ as the result of replacing every monomial of degree greater than one in $\phi(\x)$ with the corresponding RLT variable. Then, the linear relaxation of problem \ref{eq:PO} is given by
\begin{equation*}\label{eq:PPlin}\tag*{$LP(\Omega)$}
\begin{array}{rll}
	\text{minimize} & [\phi_0(\x)]_L &  \\
	\text{subject to} & [\phi_r(\x)]_L\geq \beta_r, & r=1,\ldots,R_1  \\
		         & [\phi_r(\x)]_L=\beta_r, & r=R_1+1,\ldots,R  \\
						 & [\prod_{j\in J_1}(x_j-l_j)\prod_{j\in J_2}(u_j-x_j)]_L \geq 0, &\begin{matrix}\forall J_1\cup J_2\subset N^\delta, \\ \abs{J_1\cup J_2}=\delta \end{matrix}\\
						 & \x \in \Omega \subset \mathbb{R}^{|N|}. &
\end{array}
\end{equation*}
The core idea of RLT is to embed these linear relaxations into a branch-and-bound algorithm. The constraints of the form $F_\delta(J_1,J_2)=\prod_{j\in J_1}(x_j-l_j)\prod_{j\in J_2}(u_j-x_j)\geq 0$ are called bound-factor constraints. They are introduced in order to get tighter relaxations and ensure the convergence of the RLT-based algorithm to a global optimum of~\ref{eq:PO}. \cite{Sherali1992} show that the number of bound-factor constraints in \ref{eq:PPlin} is the combinatorial number $2|N|+\delta-1$ over $\delta$. This number grows exponentially fast on $|N|$ and $\delta$ and, thus, the relaxations become computationally unmanageable as the degree and the number of variables in the problem increase. In order to mitigate this limitation, \cite{Dalkiran2013} identify a particular collection of monomials, called $J$-sets, and show that convergence is still guaranteed if only the bound-factor constraints associated with these monomials are added to the linear relaxation. This approach can dramatically reduce the size of the resulting relaxations (at the price of getting slightly looser relaxations).


A crucial step in any branch-and-bound algorithm is the selection of the branching variable. To this aim, given the optimal solution at a given node, we assign a score $\theta_j$ to each variable $j\in N$, considering the violations of the RLT identities \eqref{eq:RLTidentity} in which it is involved. Then, the variable with the highest score is selected for branching. \cite{Sherali1992} compute these scores using the following formula:
\begin{equation}
\displaystyle \theta_j=\max_{J\subset N^\delta , 1\leq \abs{J}\leq \delta -1}\Big\{\abs{\bar{X}_{J\cup\{j\}}-\bar{x}_j\bar{X}_J}\Big\}\text{, for each }j\in N,
\label{eq:branching}
\end{equation}
where $(\xb,\Xb)$ is the optimal solution of the corresponding node. More recent contributions such as \cite{Sherali2011,Sherali2011a} and \cite{GonzalezRodriguez2023} rely on variations of~\eqref{eq:branching} where the maxima are replaced with sums. In particular, \cite{GonzalezRodriguez2023} define a family of branching rules based on the formula:
\begin{equation}
\displaystyle \theta_j=\sum_{J\subset N^\delta, 1\leq \abs{J}\leq \delta -1}w(j,J)\cdot\abs{\bar{X}_{J\cup\{j\}}-\bar{x}_j\bar{X}_J}\text{, for each }j\in N,
\label{eq:newbranching}
\end{equation}
where $w(j,J) > 0$ for each $j\in N$ and $J\subset N^\delta$, with $1\leq \abs{J}\leq \delta -1$. Note that, when only  $J$-sets are used, Equations \eqref{eq:branching} and~\eqref{eq:newbranching} can be readily extended since, for each $j\in N$, we can restrict in \eqref{eq:branching} and~\eqref{eq:newbranching} to consider only the monomials $J\subset N^\delta$, with $1\leq \abs{J}\leq \delta$, such that $J\cup\{j\}$ appears in the linear relaxation built using the $J$-sets. \cite{GonzalezRodriguez2022a} proves the convergence of the RLT-based algorithm under any branching rule derived from~\eqref{eq:newbranching} and \cite {GonzalezRodriguez2023} and \cite{Ghaddar2023} present a detailed analysis of the performance of several branching rules of this family, including the following ones:
\begin{itemize}
    \item {\sc sum.} For each $j\in N$, $\theta_j$ is just the raw sum of the violations in~\eqref{eq:newbranching}, \emph{i.e.}, all $w(j,J)$ weights are set to one.
    \item {\sc range.} For each $j\in N$, $w(j,J)$ is given by the quotient between the range of the variable at the current node and its range at the root node (higher priority is given to variables whose range has been reduced less).
    \item {\sc dual.} The weights $w(j,J)$ correspond with the sum of the absolute values of the dual values associated with the constraints of the problem containing $J\cup \{j\}$.
\end{itemize}

RLT, as introduced in \cite{Sherali1992} and later implemented in \cite{Dalkiran2016}, was conceived to solve polynomial optimization problems of the form~\ref{eq:PO}, where all variables are continuous. The main contribution of this paper is to propose and a natural extension to accommodate integer variables and to analyze its performance as a function of a number of basic algorithmic decisions. We conclude this section by presenting, in Figure~\ref{fig:RLTalg}, the general scheme of the baseline (continuous) RLT algorithm. It proceeds by defining a series of linear relaxations, $LP(\Omega^k)$, of polynomial optimization problems analogous to~\ref{eq:PO}, but with respect to different hyperrectangles, $\Omega^k$. For $j \in N$ and $k \in \N$, $l^k_j$ and $u^k_j$ denote the lower and upper bound of variable $x_j$ in $LP(\Omega^k)$. 

\begin{figure}[!htbp]
\centering
\begin{tcolorbox}
\small
\textbf{Initialization.} Let $LB\vcentcolon=-\infty$, $UB\vcentcolon=+\infty$, and $\Omega^1 \vcentcolon= \Omega$. Let $LB^1 \vcentcolon= -\infty$ be the root node's potential and let $\cola\vcentcolon=\{1\}$. Let $\itercounter\vcentcolon=1$ and let $\x^{best}\vcentcolon=\emptyset$.

\textbf{Stage 1 \emph{(main)}.} Choose $k\in \cola$ such that $LB^k=\min_{s\in \cola} LB^s$. Let $\cola \vcentcolon= \cola \backslash \{ k\}$. Solve problem $LP(\Omega^k)$.
	\begin{shortitem}
		\item If $LP(\Omega^k)$ is infeasible, go to {\bf Stage~2}.
		\item If $LP(\Omega^k)$ is feasible, let $\zbar^k$ be the optimal value and let $(\xb^k, \Xb^k)$ be an optimal solution.
	\begin{shortitem}
		\item If $\zbar^k<UB$ and $\theta_j^k = 0$ for all $j \in N$:
		\begin{shortitem}
			\item \emph{Update upper bound.} $UB\vcentcolon=\zbar^k$. Let $\x^{best}\vcentcolon=\xb^k$.
			\item \emph{Prune.} Remove all $s\in \cola$ such that $LB^s \geq UB$. Go to {\bf Stage~2}.
		\end{shortitem} 
		\item If $\zbar^k < UB$ and $\theta_j^k > 0$ for some $j \in N$:
		\begin{shortitem}
			\item \emph{Branch.} Choose $p \in N$ such that $\theta_p^k = \max_{j \in N} \theta_j^k$. Branch at $\beta \vcentcolon= \bar x^k_p$, by defining $\Omega^{\itercounter+1} \vcentcolon= \Omega^k \cap \{ \x \in \R^{\vert N \vert} \text{ s.t. } l^k_p \leq x_p \leq \beta \}$ and $\Omega^{\itercounter+2} \vcentcolon= \Omega^k \cap \{ \x \in \R^{\vert N \vert} \text{ s.t. } \beta \leq x_p \leq u^k_p \}$.
	\item Update queue $\cola \vcentcolon= \cola \cup \{ \itercounter+1, \itercounter+2\}$. Let $LB^{\itercounter+1} = LB^{\itercounter+2} \vcentcolon= \zbar^k$. Let $\itercounter \vcentcolon= \itercounter + 2$. Go to {\bf Stage~2}.
		\end{shortitem}
	\item If $\zbar^j \geq  UB$. \emph{Prune branch.} Go to {\bf Stage~2}.
	\end{shortitem}
\end{shortitem}

\textbf{Stage 2 \emph{(control)}.}
\emph{Update lower bound.} $LB\vcentcolon=\min\{\min_{k\in \cola}{LB^k}, UB \}$.
\begin{shortitem}
	\item If $LB=UB$, \stop:
	\begin{shortitem}
		\item If $\x^{best}=\emptyset$, \ref{eq:PO} is infeasible.
		\item If $\x^{best}\neq\emptyset$, $\x^{best}$ is a global optimum and $UB$ is the optimal value.
	\end{shortitem} 
	\item In other case, go to {\bf Stage~1}.
\end{shortitem}
\end{tcolorbox}
\caption{Scheme of the baseline (continuous) RLT-based algorithm.}
\label{fig:RLTalg}
\end{figure}

\subsection{Computational setup}
In this section we present the computational setup in which we assess the performance of different approaches to adapt the RLT implementation in \cite{GonzalezRodriguez2023} to handle integer variables.

\subsubsection{The RLT-based solver: \solver{RAPOSa}}\label{sec:raposa}
\paragraph{Baseline solver.} The RLT-based solver employed for implementing the algorithmic ideas presented in this paper is \solver{RAPOSa} (\textbf{R}eformulation \textbf{A}lgorithm for \textbf{P}olynomial \textbf{O}ptimization - \textbf{Sa}ntiago). \solver{RAPOSa} is an RLT-based solver introduced in \cite{GonzalezRodriguez2023}. It is implemented in \CC and it is freely distributed and available for Linux, Windows, and macOS. \solver{RAPOSa} can be run from AMPL \citep{ampl}, Pyomo \citep{pyomo}, JuMP \citep{jump}, and NEOS server \citep{neos}.\footnote{For more details regarding \solver{RAPOSa} and its capabilities, please refer to \url{https://raposa.usc.es}.}

For our analysis, we build upon \solver{RAPOSa}~4.2.0, which is equipped with some well known RLT enhancements. In particular,  it includes i)~the use of $J$-sets, ii)~\solver{Gurobi} \citep{gurobi} as the auxiliary (mixed-integer) linear solver and \solver{Ipopt} \citep{Waechter2005} as the auxiliary nonlinear solver, iii)~Bound tightening \citep{Belotti2009,Puranik2017}: both OBBT, optimization-based bound tightening (at the root node) and FBBT, feasibility-based bound tightening (at all nodes), iv)~duality based branching rule, {\sc dual}, and vi)~extra tightening through products of constraint factors and bound factors. 

\paragraph{Branching point selection.} We briefly discuss now the default approach in \solver{RAPOSa}~4.2.0 regarding the selection of the branching point, since it will be relevant for the discussion of some of the results in Section~\ref{sec:bounds}. According to the baseline RLT-based algorithm in Figure~\ref{fig:RLTalg}, once subproblem~$k$ has been solved and a branching variable~$p$ has been chosen, branching is performed at the value of the variable in the optimal solution, $\bar x^k_p$.  However, it has been acknowledged by past literature that it might be beneficial to branch on other value from the range of the variable at that node $[l^k_p,u^k_p]$, a natural approach being to take some convex combination of $\bar x^k_p$ and the midpoint of the range of the variable at that node, \emph{i.e.}, $(l^k_p+u^k_p)/2$. Doing so might lead to more balanced trees and also prevent the algorithm from repeatedly branching very close to the lower or upper bound of a given variable (refer, for instance, to \cite{Belotti2009} and \cite{Ryoo:1996} for deeper discussions). Approaches along these lines are followed in the spatial branching of state-of-the-art of the art solvers such as \solver{BARON} \citep{Khajavirad2018,baron}, \solver{Couenne} \citep{Belotti2009}, \solver{Gurobi} \citep{gurobi}, and \solver{Xpress} \citep{xpress}, although the specifics, such as the precise coefficients for the convex combination, may differ across solvers. In the particular case of \solver{RAPOSa}~4.2.0, the branching point is chosen as $\beta:=0.75\cdot \bar x^k_p+0.25 \cdot (l^k_p+u^k_p)/2$, which proved to be among the best performing choices in \cite{gomezcasares2024}. In some of the above references it has also been argued that, in order to close the gap more quickly, it may be beneficial to branch on the value of the branching variable in the best available solution, $\x^{best}$, whenever possible, \emph{i.e.}, whenever $\x^{best}\neq \emptyset$ and $\x^{best}_p\in [l^k_p,u^k_p]$. Building again upon the analysis in \cite{gomezcasares2024}, this is the default behavior in \solver{RAPOSa}~4.2.0.

\paragraph{Moving to mixed-integer problems.} When moving from continuous to mixed-integer problems, an important aspect that needs to be addressed is the role of the auxiliary local solver in the RLT-based scheme. When solving continuous problems, \solver{RAPOSa} calls a continuous nonlinear local solver at selected nodes of the branch-and-bound tree, in order to get better upper bounds on the optimal value of the original problem and speed up convergence. It is worth noting that these calls give the local solver the solution of the current relaxation as starting point, so that different calls may lead to different local optima. Since we now have a mixed-integer problem, a continuous solver is not be capable of providing valid upper bounds consistently (the associated solutions can violate the integrality constraints) but, on the other hand, using instead a mixed-integer local solver is computationally more demanding. We start the analysis with a baseline version of \solver{RAPOSa} in which the solution provided by the continuous local solver is rounded and, if the resulting solution is feasible, the corresponding upper bound is computed. Later on, Section~\ref{sec:ubs} is devoted to carefully discuss different approaches to regularly obtaining upper bounds. Some of these approaches rely on calls to an auxiliary mixed-integer nonlinear local solver. In the computational analysis, \solver{Knitro} \citep{Byrd2006} is the solver used by \solver{RAPOSa} for these calls.

\subsubsection{Testing environment}
All the computational analyses reported in this paper have been performed on the supercomputer Finisterrae~III, provided by Galicia Supercomputing Centre (CESGA). More specifically, the compute nodes employed are equipped with 32 cores Intel Xeon Ice Lake 8352Y CPUs, each with 256 GB of RAM. These nodes are connected via an Infiniband HDR network, and each one is powered with 1 TB of SSD storage.

The analysis is conducted using two different sets of problems. The first one contains the mixed-integer polynomial optimization problems with box-constrained variables from MINLPLIB \citep{Bussieck2003}. The second one contains the mixed-integer quadratic optimization problems with box-constrained variables from QPLIB \citep{Furini2018}. Since we are extending RLT to optimization problems with integer variables, we select only the instances that include at least one integer variable. As a result, we have one test set from MINLPLIB and another one from QPLIB and, coincidentally, both of them have 139 instances.\footnote{There were a total of 40 problems that were disregarded because \solver{RAPOSa} was not even able to preprocess them, because of their size.}

All executions were run taking as stopping criterion that the relative or absolute gap is below the threshold $0.001$. Time limit was set to one hour.

\subsubsection{Reporting of numerical results}

The primary reporting tool consists of a series of summary tables, divided into two subtables each (one for MINLPLIB and one for QPLIB), devoted to assess the performance of different approaches according to the following metrics:

\begin{description}
\item[Unsolved.] Number of unsolved instances. In brackets we show the number of instances not solved by any approach, followed by the total number of instances in the corresponding test set.
\item[$\boldsymbol{\text{Gap} = \infty}$.] Number of instances in which the algorithm terminated without providing an optimality gap. In brackets we show number of instances for which no approach delivered an optimality gap.
\item[$\boldsymbol{\text{UB} = \infty}$.] Number of instances in which the algorithm terminated without providing an upper bound. In brackets we show number of instances for which no approach delivered an upper bound.
\item[$\boldsymbol{\text{LB} = -\infty}$.] Number of instances in which the algorithm terminated without providing a lower bound. In brackets we show number of instances for which no approach delivered a lower bound.
\item[Gap.] Geometric mean of the relative optimality gap.\footnote{The adoption of geometric means is becoming the standard in benchmarking, due to their lower sensitivity to outliers compared to the standard mean. Note, however, that this choice tends to yield more condensed results, potentially diminishing the perceived impact of improvements.} We exclude instances solved by all approaches and instances for which at least one approach could not return a gap after the time limit. In brackets we show the remaining number of instances.
\item[Time.] Geometric mean of the running time, but excluding those instances solved by all approaches in less than 5~seconds and also those not solved by any approach within the time limit. In brackets we show the remaining number of instances.
\item[Nodes.] Geometric mean of the number of nodes generated by the RLT-based algorithm in those instances solved by all approaches. In brackets we show the remaining number of instances.
\end{description}

Note that, for all the above metrics, the smaller the value, the better. The first column in all tables represents a baseline or reference configuration on which to assess the relative performance of the alternative approaches. To facilitate the interpretation of the results, when reporting results for ``Gap'', ``Time'', and ``Nodes'', we use absolute values in the first column and, in the rest, percentages representing variations with respect to the baseline. On the other hand, since ``Unsolved'', ``$\text{Gap} = \infty$'', ``$\text{UB} = \infty$'', and ``$\text{LB} = -\infty$'' refer to numbers of instances, we report them with absolute values in all columns.

\section{RLT in mixed-integer polynomial optimization}\label{sec:RLTint}

The primary objective of this section is to present a first functional and robust extension of RLT, as introduced in \cite{Sherali1992}, to mixed-integer polynomial optimization problems. Despite the abundant literature around RLT, both for the continuous and the discrete setting, such an extension has been minimally explored in the literature. One exception is \cite{GonzalezRodriguez2022a}, where it is proven that the same RLT-based algorithm, but solving MILP relaxations instead of just LP ones, converges to a global optimum of the mixed-integer polynomial optimization problem. We delve into a detailed discussion of this approach in Section~\ref{milpsolver}. Before doing so, in Section~\ref{sec:intbranch} we discuss alternative modifications of the baseline algorithm in Figure~\ref{fig:RLTalg}.

It is important to note that, in all the approaches discussed in this work, continuous and integer variables are treated equally in the reformulation phase. This implies that both types of variables are used in the definition of the bound-factor constraints, which allows the tightness of the latter to be fully exploited. On the downside, given that the number of bound-factor constraints grows exponentially fast on the number of variables and the degree of the monomials, the above approach might lead to exceedingly large subproblems for some instances.

An alternative approach, which is beyond the scope of this paper, would be to build upon alternative reformulations of the discrete variables, along the lines of the different references mentioned in the introduction. In particular, many of these references discuss linearizations of products of binary variables, the introduction of valid inequalities, and the tightness of the resulting relaxations. Even if such relaxations are looser than the ones provided by the approaches studied here, they might lead to a significantly smaller number of bound-factor constraints, which might in turn improve performance. Finding the right balance between the two types of reformulations is far from trivial and is definitely a must for future research, especially given the prevalence of products of binary variables in many real-world applications.

\subsection{Prioritization of branching decisions}\label{sec:intbranch}

Consider the RLT-based branch-and-bound algorithm for polynomial optimization problems outlined in Section~\ref{sec:RLT}, and whose general scheme is depicted in Figure~\ref{fig:RLTalg}. Whenever we solve a linear relaxation at a specific node, the only violations preventing the solution of that node from being feasible in the original polynomial optimization problem are the RLT violations defined in~\eqref{eq:RLTidentity}. In the presence of integer variables, the violations of the integrality constraints have to be taken into account as well. Thus, when defining a branching rule we have to decide how to prioritize the violations of the RLT identities and of the integrality constraints. As noted in the introduction, this prioritization is handled differently by different modern solvers. For instance, \solver{Couenne} \citep{Belotti2009}, in its default configuration, gives full priority to branching on integer variables and \solver{Xpress} \citep{xpress} determines the branching variable by computing scores along the lines of reliability branching \citep{Achterberg2005}, with no particular priority given to either continuous or integer variables. 

There are many possibilities to jointly handle both types of violations and, yet, despite the huge impact these choices may have in the performance of the algorithm, we have found no reference studying and comparing them. In this section we concentrate on two of extreme choices. The first one, RLT-first, involves prioritizing the violations of the RLT identities. Only when the RLT identities at a given node are zero, do we consider the violations of the integrality constraints. The second approach, integrality-first, does the opposite: it uses the information regarding the violations of the integrality constraints first, and only when the integrality violations are zero, does it consider the RLT violations. In order to prove the convergence of the resulting branch-and-bound algorithms to the global optimum of the original problem one can easily adapt the arguments of Theorem~1.9 in \cite{GonzalezRodriguez2022a}. 

For this preliminary analysis comparing RLT-first and integrality-first, when considering integrality violations, we just branch on the variable that is furthest away from being integer. 

\begin{table}[htbp]
\centering
{\begin{tabular}{l|lrr}
\toprule
     Test set &               &  Integrality-first & RLT-first\\
\midrule
					& Unsolved \nf{(61/139)} &  $64\phantom{.00}$& $	70\phantom{.00}$\phantom{\%} \\
					& $\text{Gap} = \infty$ \nf{(38)} &        $39\phantom{.00}$   &  $50\phantom{.00}$\phantom{\%} \\
					& $\text{UB} = \infty$ \nf{(38)} &              $39\phantom{.00}$       &   $50\phantom{.00}$\phantom{\%}  \\
MINLPLIB	& $\text{LB} = -\infty$ \nf{(0)} &         $0\phantom{.00}$       &   $0\phantom{.00}$\phantom{\%}\Bstrut \\ \cdashline{2-4}
					& Gap \nf{(22)} &   $0.30$   &  $-32.58$\%\Tstrut \\
					& Time \nf{(33)} &    $83.66$& $+47.56$\% \\
					& Nodes \nf{(66)} &   $266.81$& $-40.13$\% \\
\midrule
					& Unsolved \nf{(130/139)} &  $134\phantom{.00}$& $	132\phantom{.00}$\phantom{\%} \\
					& $\text{Gap} = \infty$ \nf{(34)} &        $37\phantom{.00}$   &  $37\phantom{.00}$\phantom{\%} \\
					& $\text{UB} = \infty$ \nf{(34)} &              $37\phantom{.00}$       &   $37\phantom{.00}$\phantom{\%}  \\
QPLIB			& $\text{LB} = -\infty$ \nf{(0)} &         $0\phantom{.00}$       &   $0\phantom{.00}$\phantom{\%}\Bstrut \\ \cdashline{2-4}
					& Gap \nf{(96)} &      $3.35$   &  $-19.76$\%\Tstrut \\
					& Time \nf{(9)} &    $1882.20$& $-30.18 $\% \\
					& Nodes \nf{(3)} &   $1586.66$& $-29.26$\% \\
\bottomrule
\end{tabular}}
\caption{Prioritization of integrality
 violations versus prioritization on RLT violations.}\label{table:branching}
\end{table}

The results in Table~\ref{table:branching} show that none of the approaches dominates the other. The unsolved instances are similar for the two approaches, with integrality-first performing slightly better. Regarding gaps, we can see that, in both test sets, in the instances where both approaches ended up the running time with a finite gap, RLT-first is better at closing the gap. However, integrality-first is able to deliver a finite gap in 11 instances more than RLT-first. In terms of running times, integrality-first does better in MINLPLIB and RLT-first does better in QPLIB. Yet, it is worth noting that in the later test set these times are computed over just 9 instances. Finally, the number of nodes suggests that RLT-first is better at reducing the size of the branch-and-bound tree.

Even though there is not a clear winner between RLT-first and integrality-first, the latter has an important advantage in terms of implementation: it can rely on a state-of-the-art MILP solver to handle the relaxations. Thus, instead of solving linear relaxations of the form \ref{eq:PPlin} at each iteration, we can solve the resulting mixed-integer linear relaxation by adding the integrality constraints to \ref{eq:PPlin}. With respect to the RLT scheme in Figure~\ref{fig:RLTalg}, the only change is precisely that the relaxations $LP(\Omega^k)$ are now MILP problems. This allows the RLT-based algorithm to benefit from all the capabilities of a state-of-the-art MILP solver in terms of MILP reformulations, cuts, branching decisions,\ldots In Table~\ref{tab:branch_vs_milpsolver} we compare the performance of RLT-first and integrality-first with this new version in which integrality-first is implemented by letting \solver{Gurobi} \citep{gurobi} solve the MILP relaxations starting already at the root node (depth $d=1$).\footnote{We do not include Nodes in Table~\ref{tab:branch_vs_milpsolver}, since option ``$d=1$'' will lead to much smaller trees in \solver{RAPOSa} just because all the integer branching is handled by the MILP solver.} The superiority of the new version is outstanding. The only weakness is in the number of instances in which it was not able to find a lower bound, just because the root-node MILP relaxation was not solved within the time limit. As we discuss in Section~\ref{sec:implbs}, this can be easily addressed by first solving the LP relaxation at the root node to compute an initial lower bound.

\begin{table}[!htpb]
\centering
{\begin{tabular}{l|lrrr}
\toprule
     Test set &               &  \multicolumn{1}{c}{$d=1$} &  Integrality-first & RLT-first\\
\midrule
					& Unsolved \nf{(26/139)} & 26\phantom{.000} & 64\phantom{.00}\phantom{\%}&	70\phantom{.00}\phantom{\%}	\\
					&$\text{Gap} = \infty$ \nf{(16)} & 20\phantom{.000} & 39\phantom{.00}\phantom{\%}	&50\phantom{.00}\phantom{\%} \\
					&$\text{UB} = \infty$ \nf{(16)} & 20\phantom{.000} & 39\phantom{.00}\phantom{\%}	&50\phantom{.00}\phantom{\%} \\
MINLPLIB	&$\text{LB} = -\infty$ \nf{(0)} & 6\phantom{.000} & 0\phantom{.00}\phantom{\%}& 0\phantom{.00}\phantom{\%}\Bstrut \\ \cdashline{2-5}
					& Gap \nf{(19)} & 0.006 & +4390.95\% & +2756.67\%\Tstrut \\
					& Time \nf{(67)} & 11.661 & +4707.75\% & +5762.53\% \\
	\midrule
					& Unsolved \nf{(63/139)} & 65\phantom{.000} & 134\phantom{.00}\phantom{\%}&	132\phantom{.00}\phantom{\%}	\\
					&$\text{Gap} = \infty$ \nf{(20)} & 64\phantom{.000} & 37\phantom{.00}\phantom{\%}&	37\phantom{.00}\phantom{\%}	 \\
					&$\text{UB} = \infty$ \nf{(20)} & 26\phantom{.000} & 37\phantom{.00}\phantom{\%}	&37\phantom{.00}\phantom{\%} \\
QPLIB			&$\text{LB} = -\infty$ \nf{(0)} & 52\phantom{.000} & 0\phantom{.00}\phantom{\%}	&0\phantom{.00}\phantom{\%}\Bstrut \\ \cdashline{2-5}
					& Gap \nf{(55)} & 0.001 & $>10^5\,$\%\phantom{.00} & $>10^5\,$\%\phantom{.00}\Tstrut\\
					& Time \nf{(76)} & 250.664& +1222.82\% & +1174.63\%  \\
\bottomrule
\end{tabular}}
\caption{Letting a state-of-the-art MILP solver handle the integer branching.}\label{tab:branch_vs_milpsolver}
\end{table}

It is not surprising that this new version outperforms both integrality-first and RLT-first, since these were just two simple implementations intended to provide some initial results on the performance of the two opposite approaches regarding branching prioritization. Yet, the difference in performance is so big that, hereafter, we just consider RLT-based approaches in which integer branching is handled by the auxiliary MILP solver.


\subsection{Controlling integer branching depth}\label{milpsolver}


The purpose of this section is to explore whether it might be beneficial to start the algorithm with RLT-first until some depth is reached and then switch to integrality-first (with Gurobi solving the MILP subproblems). The idea is that RLT violations might be informative enough to guide the construction of the branch-and-bound tree in a first phase of the algorithm and then switch to the computationally more demanding solution of the MILP subproblems. Thus, for a given depth $d$, we solve the linear relaxation for the nodes with depth smaller than $d$ and solve the mixed-integer linear relaxation in the nodes with depth larger than or equal to $d$. 

As far as spatial branching is concerned, the branching variable is selected as usual in the RLT-based algorithm, \emph{i.e.}, scores are computed following~\eqref{eq:newbranching} and we can then follow one of the rules described in Section~\ref{sec:prelim}: {\sc sum}, {\sc range}, or {\sc dual}. There is an important aspect we need to highlight. In the analysis in \cite{Ghaddar2023}, the best performing branching rule was {\sc dual}, which computes the $w(j, J)$ weights as the sum of the absolute values of the dual values associated with the constraints of the problem containing $J\cup \{j\}$. Nevertheless, when we let a MILP solver handle the MILP relaxations, this dual values are not available anymore. In order to circumvent this problem, in this section we switch to {\sc sum} as soon as the MILP solver enters. We discuss this in more detail in Section~\ref{sec:dual}, where we explore some alternatives.

\begin{table}[!htpb]
\centering
{\begin{tabular}{l|lrrrr}
\toprule
     Test set &               &  \multicolumn{1}{c}{$d=1$} & \multicolumn{1}{c}{$d=5$}& \multicolumn{1}{c}{$d=10$}& \multicolumn{1}{c}{$d=50$}  \\
\midrule
        & Unsolved \nf{(26/139)} & 26\phantom{.000} & 37\phantom{.00}\phantom{\%}&	42\phantom{.00}\phantom{\%}	&70\phantom{.00}\phantom{\%}\\
        &$\text{Gap} = \infty$ \nf{(15)} & 19\phantom{.000} & 28\phantom{.00}\phantom{\%}	&27\phantom{.00}\phantom{\%}&	49\phantom{.00}\phantom{\%} \\
        &$\text{UB} = \infty$ \nf{(15)} & 19\phantom{.000} & 28\phantom{.00}\phantom{\%}	&27\phantom{.00}\phantom{\%}	&49\phantom{.00}\phantom{\%} \\
MINLPLIB&$\text{LB} = -\infty$ \nf{(0)} & 6\phantom{.000} & 0\phantom{.00}\phantom{\%}& 0\phantom{.00}\phantom{\%} & 0\phantom{.00}\phantom{\%}\Bstrut \\ \cdashline{2-6}
				& Gap \nf{(18)} & 0.006 & +43.84\% & +431.34\% & +4242.67\%\Tstrut \\
        & Time \nf{(62)} & 14.565 & +280.79\% & +976.25\% & +7534.46\%\\
        & Nodes \nf{(69)} & 6.964 & +431.16\% & +1559.96\% & +2869.63\% \\
\midrule
	       & Unsolved \nf{(63/139)} & 65\phantom{.000} & 76\phantom{.00}\phantom{\%}&	102\phantom{.00}\phantom{\%}	&132\phantom{.00}\phantom{\%} \\	
        &$\text{Gap} = \infty$ \nf{(20)} & 64\phantom{.000} & 24\phantom{.00}\phantom{\%}&	25\phantom{.00}\phantom{\%}	&37\phantom{.00}\phantom{\%} \\
        &$\text{UB} = \infty$ \nf{(20)} & 26\phantom{.000} & 24\phantom{.00}\phantom{\%}	&25\phantom{.00}\phantom{\%}&	37\phantom{.00}\phantom{\%} \\
QPLIB   &$\text{LB} = -\infty$ \nf{(0)} & 52\phantom{.000} & 0\phantom{.00}\phantom{\%}	&0\phantom{.00}\phantom{\%}&	0\phantom{.00}\phantom{\%}\Bstrut \\ \cdashline{2-6}
        & Gap \nf{(53)} & 0.001 & +463.47\% & +17878.58\% & $>10^5\,$\%\phantom{.00}\Tstrut\\
				& Time \nf{(76)} & 250.664& +159.29\% & +638.60\% & +1174.65\% \\
        & Nodes \nf{(7)} & 1\phantom{.000} & +2334.05\% & +38191.25\% & $>10^5\,$\%\phantom{.00} \\
\bottomrule
\end{tabular}}
\caption{Depth from which a MILP solver starts to solve MILP relaxations.}\label{table:gurobi}
\end{table}

Table~\ref{table:gurobi} shows a comparison for different values of the depth $d$ at which the MILP solver starts to solve the MILP relaxations. Performance clearly deteriorates as we increase the depth at which the MILP solver enters. Therefore, the approach with the best performance is to let the MILP solver handle the MILP relaxations from the beginning of the branch-and-bound tree. It is worth noting that in QPLIB, the 7 instances solved by all the approaches were solved at the root node when we solve directly the MILP relaxation, \emph{i.e.}, $d=1$.
Again, the main weakness for the approach $d=1$ comes from metric $LB=-\infty$, which, as mentioned when discussing Table~\ref{tab:branch_vs_milpsolver}, can be addressed by first solving the LP relaxation at the root node to compute an initial lower bound (Section~\ref{sec:implbs}).

\subsection{Comparison with other state-of-the-art solvers}\label{sec:othersolvers}

To conclude this section, we compare the performance of the mixed-integer RLT-based algorithm with two state-of-the-art MINLP solvers: \solver{BARON} \citep{baron} and \solver{Couenne} \citep{Belotti2009}. The main purpose of this section is to assess the strengths and weaknesses of our extension of RLT to mixed-integer polynomial optimization problems. In Table \ref{table:solvers1}, we present a comparison between \solver{RAPOSa}~4.2.0 ($d=1$), which is our implementation of the RLT-based algorithm in \solver{RAPOSa}, where the MILP solver starts to solve MILP relaxation at the root node, \solver{BARON} 23.3.1, and \solver{Couenne} 0.5.7. The auxiliary MILP solver used by \solver{RAPOSa} is \solver{Gurobi} \citep{gurobi}. Regarding \solver{BARON}, we include two different versions: the default one, which uses \solver{CBC} \citep{JohnForrest2023} as the auxiliary MILP solver, and another one that uses \solver{CPLEX} \citep{cplex}. This last version is important to have a fair comparison with \solver{RAPOSa} using \solver{Gurobi}. For \solver{Couenne}, we include in our analysis the default version that uses \solver{CBC}. Differently from the results in the previous subsections, where we imposed no limit on the number of threads used by \solver{RAPOSa}, we now impose a common limit of 8 threads to all solvers. This choice leaves \solver{Couenne} in a disadvantageous position, since it does not have parallelization capabilities, whereas \solver{RAPOSa} and \solver{BARON} lean on the parallelizations behind the auxiliary MILP solvers (neither \solver{RAPOSa} nor \solver{BARON} perform parallelization in the spatial branching).

\begin{table}[!htpb]
\centering
{\begin{tabular}{l|lrrrr}
\toprule
     Test set &          &  \solver{RAPOSa} &  \solver{Couenne} & \solver{BARON} (\solver{CBC})& \solver{BARON} (\solver{CPLEX})  \\
\midrule
					& Unsolved \nf{(17/139)} & 30\phantom{.000}& 41\phantom{.00}\phantom{\%} &30\phantom{.00}\phantom{\%}&	25\phantom{.00}\phantom{\%} \\
					&$\text{Gap} = \infty$ \nf{(0)} & 22\phantom{.000} & 1\phantom{.00}\phantom{\%} & 0\phantom{.00}\phantom{\%} & 1\phantom{.00}\phantom{\%} \\
					&$\text{UB} = \infty$ \nf{(0)} & 22\phantom{.000} & 1\phantom{.00}\phantom{\%}& 0\phantom{.00}\phantom{\%} & 1\phantom{.00}\phantom{\%}  \\
MINLPLIB	&$\text{LB} = -\infty$ \nf{(0)} & 7\phantom{.000} & 0\phantom{.00}\phantom{\%}& 0\phantom{.00}\phantom{\%} & 0\phantom{.00}\phantom{\%}\Bstrut \\ \cdashline{2-6}
					& Gap \nf{(27)} & 0.007& +3256.06\%& +276.21\% & +83.90\%\Tstrut  \\
					& Time \nf{(67)} & 29.794& +503.68\% & +35.83\% & $-77.33\%$\\
					\midrule
					& Unsolved \nf{(55/139)} & 67\phantom{.000} & 131\phantom{.00}\phantom{\%}& 105\phantom{.00}\phantom{\%}	& 79\phantom{.00}\phantom{\%}	 \\
					&$\text{Gap} = \infty$ \nf{(0)} & 66\phantom{.000}&6\phantom{.00}\phantom{\%} &0\phantom{.00}\phantom{\%}&	0\phantom{.00}\phantom{\%}	 \\
					&$\text{UB} = \infty$ \nf{(0)} & 27\phantom{.000}&	6\phantom{.00}\phantom{\%} &0\phantom{.00}\phantom{\%}	&0\phantom{.00}\phantom{\%} \\
QPLIB			&$\text{LB} = -\infty$ \nf{(0)} & 55\phantom{.000} & 1\phantom{.00}\phantom{\%}& 0\phantom{.00}\phantom{\%} & 0\phantom{.00}\phantom{\%}\Bstrut \\ \cdashline{2-6}
					& Gap \nf{(69)} & 0.001 & +41269.26\% & +3467.69\% & +895.42\%\Tstrut\\
					& Time \nf{(84)} & 396.682& +677.72\% & +241.57\% & $-22.72\%$  \\
\bottomrule
\end{tabular}}
\caption{Comparison of \solver{RAPOSa} with other state-of-the-art solvers.}\label{table:solvers1}
\end{table}

Table~\ref{table:solvers1} presents a comparison between the different solvers in both test sets. Note that we have removed the row with the number of nodes, since they are not comparable across different solvers. Regarding Unsolved, \solver{BARON} (\solver{CPLEX}) and \solver{RAPOSa} are the best performing solvers in MINLPLIB and QPLIB, respectively. If we look at Gap, \solver{RAPOSa} is by far the solver that closes the gap the most within the time limit (in those instances not solved by all solvers and for which all solvers returned a finite optimality gap). Regarding Time, the best-performing solver is the version of \solver{BARON} (\solver{CPLEX}) as an auxiliary MILP solver, followed by \solver{RAPOSa}. Regarding \solver{RAPOSa} and \solver{BARON} (\solver{CPLEX}), the numbers in Table~\ref{table:solvers1} may suggest that the difference in Gap in favor of \solver{RAPOSa} is much larger than the difference in Time in favor of \solver{BARON} (\solver{CPLEX}). Yet, this is not the case since, for instance, $+900\%$ implies a factor of~10 increase and the corresponding  factor of~10 decrease would be given by $-90\%$. Then, if we look at Gap, for MINLPLIB, the $+83.90\%$ increase means moving from a geometric mean of $0.007$ to $0.013$ and, for QPLIB, the $+895.42\%$ increase means moving from $0.0012$ to $0.012$, a factor of~10. On the other hand, regarding Time, for MINLPLIB, the $-77.33\%$ decrease means moving from $29.79$ to $6.75$, close to a factor of~5, whereas, for QPLIB, $-22.72\%$ means moving from $396.68$ to $306.56$. Thus, the superiority of \solver{RAPOSa} in Gap is comparable to the superiority of \solver{BARON} (\solver{CPLEX}) in Time.

In a nutshell, Table~\ref{table:solvers1} shows that \solver{RAPOSa} seems to be better at closing the gap in difficult instances and \solver{BARON} (\solver{CPLEX}) seems to be better at solving the not-so-difficult ones. These results are very promising, since they have been obtained with a preliminary implementation of a naive adaptation of the baseline continuous RLT-based algorithm to mixed-integer polynomial optimization problems. In recent literature, \cite{Dalkiran2016} and \cite{GonzalezRodriguez2023} showed that RLT is competitive with other techniques for continuous polynomial optimization. The results in Table~\ref{table:solvers1} show the great potential of RLT-based algorithms also for mixed-integer polynomial optimization.

Table~\ref{table:solvers1} also shows that the current version in \solver{RAPOSa} struggles to find optimality gaps. There are 88 instances out of 278 in both test sets where \solver{RAPOSa} fails to return an optimality gap within the time limit, while \solver{BARON} (\solver{CPLEX}) manages to return a gap in all instances except 1 and \solver{Couenne} fails in only 7 instances. Looking at the rows regarding lower and upper bound, we find that \solver{RAPOSa} fails to find both lower and upper bounds. Therefore, in the remainder of this paper, we focus on improving this aspect, aiming to effectively compute more lower and upper bounds without compromising the overall performance of \solver{RAPOSa}.

\begin{figure}[!htpb]
\centering
\begin{subfigure}{0.48\textwidth}
\includesvg[pretex=\scriptsize,width=\textwidth]{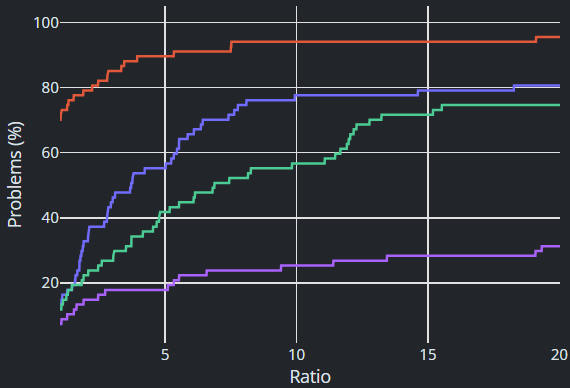}
    \caption{Running time in MINLPLIB.}
    \label{fig:timeminlpfinal}
\end{subfigure}
\hfill
\begin{subfigure}{0.48\textwidth}
    \begin{tikzpicture}
     \node (fig1) at (0,0){\includesvg[pretex=\scriptsize,width=\textwidth]{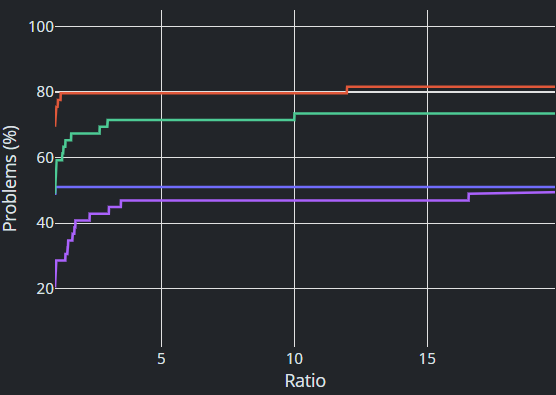}};
     \node (fig2) at (3.1,-0.9){\includesvg[pretex=\scriptsize,width=0.5\textwidth]{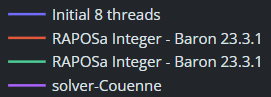}}; 
    \end{tikzpicture}
    \caption{Optimality gap in MINLPLIB.}
    \label{fig:gapminlpinit}
\end{subfigure}

\begin{subfigure}{0.48\textwidth}
    \includesvg[pretex=\scriptsize,width=\textwidth]{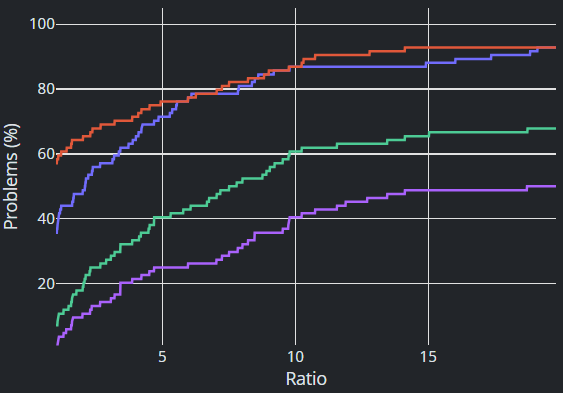}
    \caption{Running time in QPLIB.}
    \label{fig:timeqplibinit}
\end{subfigure}
\hfill
\begin{subfigure}{0.48\textwidth}
    \includesvg[pretex=\scriptsize,width=\textwidth]{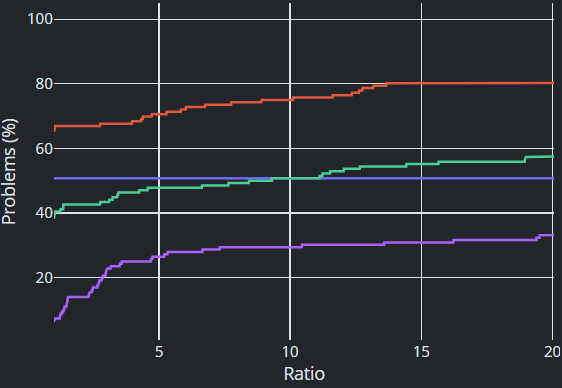}
    \caption{Optimality gap in QPLIB.}
    \label{fig:gapqplibinit}
\end{subfigure}
        
\caption{Performance profiles comparing \solver{RAPOSa} with other \solver{BARON} and \solver{Couenne}.}
\label{fig:solvers1}
\end{figure}

To conclude this section we compare the different solvers using performance profiles \citep{Dolan2002}. Figure~\ref{fig:solvers1} shows performance profiles for both running times and optimality gaps in both test sets. In the $x$-axis we represent ratios of running times or relative optimality gaps, while in the $y$-axis we represent the percentage of instances in which the corresponding configuration has a ratio lower than the value on the $x$-axis. For each instance, the ratios are computed dividing running times or relative optimality gaps of each configuration by the best configuration in that instance.\footnote{It is worth noting that the reliability of performance profiles is sometimes limited when comparing more than two configurations/solvers \citep{Gould2016}. Yet, we think that, when complemented with the numeric results in the tables, they help to obtain a more complete picture.}

The results in Figure~\ref{fig:solvers1} are consistent with those in Table~\ref{table:solvers1}. \solver{BARON} (\solver{CPLEX}) is clearly superior in Time, and particularly so in MINLPLIB. Regarding gaps, it may come as a surprise that \solver{BARON} (\solver{CPLEX}) is clearly on top of \solver{RAPOSa}, but the reason is that the latter gets heavily penalized by the large number of instances in which it is not able to return a gap after the time limit. Finally, \solver{Couenne} seems to fall behind according to both Time and Gap.

\section{Getting bounds consistently}\label{sec:bounds}

The objective of this section is to discuss various approaches to reduce the number of instances in which the RLT-based algorithm does not find lower or upper bounds. As discussed in Section \ref{sec:othersolvers}, this is crucial to address the main weakness of \solver{RAPOSa} for MINLP problems.

\subsection{Bound tightening}

In this section, we briefly discuss bound tightening techniques, which are at the core of most global optimization algorithms for nonlinear problems \citep{Belotti2009,Belotti2012a,Puranik2017}. Specifically, \cite{GonzalezRodriguez2023} and \cite{gomezcasares2024} explore the impact of optimality-based bound tightening techniques (OBBT) and feasibility-based bound tightening techniques (FBBT) on the performance of an RLT-based algorithm, demonstrating significant improvements.\footnote{Interestingly, \cite{gomezcasares2024} even explore, by building upon the analysis in \cite{raposaconic}, the possibility of performing OBBT on conic tightenings of the baseline RLT linear relaxations. More specifically, they add SOCP and SDP constraints and resort to a specialized SDP solver at the OBBT stage.} The RLT-based algorithm discussed in the previous sections ignores integrality when running OBBT and FBBT and, yet, both of them might deliver additional tightness using integrality. 

OBBT aims to obtain tighter bounds by solving a series of relaxations of minor variations of the original problem. The classic approach, implemented in \solver{RAPOSa}, involves solving a series of linear relaxations similar to \ref{eq:PPlin}, but with the objective function replaced by the minimization or maximization of each variable in the problem. In the context of MINLP problems, one may impose integrality on these linear relaxations. This introduces a trade-off between obtaining better bounds and consuming more time to solve the relaxations, as they become MILP problems instead of LP ones. In this section we study this trade-off. Recall that, as proposed in \cite{GonzalezRodriguez2023}, OBBT is performed only before solving the root node, since it is considerable time consuming. FBBT involves deriving tighter bounds directly by using propagation on the problem constraints, and is computationally much less demanding than OBBT. Integrality may be incorporated into FBBT by appropriately rounding the intermediate bounds for the integer variables during the FBBT propagation. Yet, we believe such an approach has limited potential for improvement and we have not tested it.

\begin{table}[!htpb]
\centering
{\begin{tabular}{l|lrr}
\toprule
     Test set &           &  LP OBBT & MILP OBBT \\
\midrule
						& Unsolved \nf{(26/139)} & 26\phantom{.00} & 26\phantom{.00}\phantom{\%} \\
						&$\text{Gap} = \infty$ \nf{(18)} & 19\phantom{.00} & 18\phantom{.00}\phantom{\%} \\
						&$\text{UB} = \infty$ \nf{(18)} & 19\phantom{.00} & 18\phantom{.00}\phantom{\%} \\
MINLPLIB		&$\text{LB} = -\infty$ \nf{(6)} & 6\phantom{.00} & 6\phantom{.00}\phantom{\%}\Bstrut \\ \cdashline{2-4}
						& Gap \nf{(7)} & 1.05 & +46.98\%\Tstrut \\
						& Time \nf{(37)} & 72.85 & +34.95\% \\
						& Nodes \nf{(113)} & 10.99 & $-1.61\%$ \\
						&OBBT time \nf{(139)} & 0.99 & +111.91\% \\
						\midrule
						& Unsolved \nf{(65/139)} & 65\phantom{.00} & 65\phantom{.00}\phantom{\%} \\
						&$\text{Gap} = \infty$ \nf{(64)} & 64\phantom{.00} & 64\phantom{.00}\phantom{\%} \\
						&$\text{UB} = \infty$ \nf{(26)} & 26\phantom{.00} & 26\phantom{.00}\phantom{\%} \\
QPLIB				&$\text{LB} = -\infty$ \nf{(52)} & 52\phantom{.00} & 52\phantom{.00}\phantom{\%}\Bstrut \\ \cdashline{2-4}
						& Gap \nf{(1)} & 57.02 & +0.21\%\Tstrut \\
						& Time \nf{(74)} & 231.75 & +76.88\% \\
						& Nodes \nf{(74)} & 1.24 & $-0.18\%$ \\
						& OBBT time \nf{(139)} & 70.06 & +135.23\% \\		
\bottomrule
\end{tabular}}
\caption{LP OBBT versus MILP OBBT.}\label{table:obbt}
\end{table}

In Table~\ref{table:obbt}, we present a comparison between performing OBBT solving LP relaxations and solving the MILP relaxations. First, note that we add a new row called ``OBBT time'' which is the geometric mean of the running time consumed by OBBT. If we look at Gap and Time, the performance of LP OBBT is clearly superior in both test sets. Moreover, the number of lower and upper bounds found by the two approaches is almost the same, as is the number of unsolved instances. We can explain this superior performance of LP OBBT looking at Nodes and OBBT time. First, if we look at Nodes, we see a small reduction favorable to MILP OBBT, as one would expect. This is because MILP OBBT leads to better bounds, since it solves tighter relaxations. However, if we look at the time consumed during the OBBT phase, we see that MILP OBBT consumes much more time than LP OBBT, which makes the trade-off clearly go in favor of LP OBBT. It is worth noting that we performed additional computational experiments with variations of MILP OBBT in which the subproblems were solved relaxing the tolerance imposed on the optimality gap or with smaller time limits. Although some of these variations delivered slight improvements on some metrics with respect to the version reported in Table~\ref{table:obbt}, none of them was competitive with LP OBBT. For this reason, in the following sections we stick to LP OBBT.

\subsection{Getting lower bounds}\label{sec:implbs}
We now focus on improving the RLT-based algorithm to obtain lower bounds more consistently. As discussed in Section~\ref{milpsolver}, there are instances where the RLT-based algorithm fails to return a lower bound within the running time because it cannot solve the MILP relaxation at the root node. A simple approach to obtain a lower bound in these instances would be to first solve, at the root node, the LP relaxation~\ref{eq:PPlin} before solving the MILP relaxation. While the LP relaxation will provide a lower bound equal to or worse than the one of the MILP relaxation, it is still useful to quickly obtain a lower bound.

Table~\ref{table:contrelax} presents a comparison between directly solving the MILP relaxation at the root node and solving the LP relaxation first. The performance in both test sets is very similar according most metrics: same values for Unsolved and Nodes and almost identical values for Gap and Time. Regarding Gap and Time, when solving the LP relaxation first, some time is invested in this process, leading to slightly larger running times and optimality gaps for this approach. The crucial difference, however, is that when the LP relaxation is solved at the beginning, we obtain more lower bounds in both test sets: 6 more in MINLPLIB and 52 more in QPLIB, returning  a lower bound for all instances in both test sets. Therefore, solving the LP relaxation at the beginning of the RLT-based algorithm yields better overall performance. Hereafter we adopt this approach.

\begin{table}[!htpb]
\centering
{\begin{tabular}{l|lrr}
\toprule
     Test set &           & No LP relaxation & LP relaxation \\
\midrule
						& Unsolved \nf{(26/139)} & 26\phantom{.00} & 26\phantom{.00}\phantom{\%} \\
						&$\text{Gap} = \infty$ \nf{(19)} & 19\phantom{.00} & 19\phantom{.00}\phantom{\%} \\
						&$\text{UB} = \infty$ \nf{(19)} & 19\phantom{.00} & 19\phantom{.00}\phantom{\%} \\
MINLPLIB		&$\text{LB} = -\infty$ \nf{(0)} & 6\phantom{.00} & 0\phantom{.00}\phantom{\%}\Bstrut \\ \cdashline{2-4}
						& Gap \nf{(7)} & 1.05 &  +0.00\%\Tstrut\\
						& Time \nf{(35)} & 85.05 & +0.35\% \\
						& Nodes \nf{(113)} & 10.99 &  +0.00\%\\
						\midrule
						& Unsolved \nf{(65/139)} & 65\phantom{.00} & 65\phantom{.00}\phantom{\%} \\
						&$\text{Gap} = \infty$ \nf{(26)} & 64\phantom{.00} & 26\phantom{.00}\phantom{\%} \\
						&$\text{UB} = \infty$ \nf{(26)} & 26\phantom{.00} & 26\phantom{.00}\phantom{\%} \\
QPLIB				&$\text{LB} = -\infty$ \nf{(0)} & 52\phantom{.00} & 0\phantom{.00}\phantom{\%}\Bstrut \\ \cdashline{2-4}
						& Gap \nf{(1)} & 57.02 & +0.25\%\Tstrut \\
						& Time \nf{(74)} & 231.75 & +0.44\% \\
						& Nodes \nf{(74)} & 1.24 &	+0.00\%  \\
\bottomrule
\end{tabular}}
\caption{Impact of solving the LP relaxation at the root node.}
\label{table:contrelax}
\end{table}

\subsection{Improving the branching rule}\label{sec:dual}

As we mentioned in Section~\ref{milpsolver}, once we solve MILP subproblems at the nodes of the branch-and-bound tree, we do not have the information required by the {\sc dual} branching rule, which was the best performing one for continuous problems in \cite{Ghaddar2023}. Thus, so far we have considered {\sc sum} as the branching rule, \emph{i.e.}, setting $w(j, J) = 1$ for each $j$ and $J$ in~\eqref{eq:newbranching}. In this section we study two alternative rules: an adaptation of {\sc dual} to the current setting and {\sc range}, the second-best rule in \cite{Ghaddar2023} and the best performing one in \cite{GonzalezRodriguez2023}. Motivated by what we have just seen in Section~\ref{sec:implbs}, \emph{i.e.}, the fact that solving the LP relaxation at the root node has a minimal impact on Time and Gap, we present the following adaptation of the {\sc dual} rule. At each iteration we solve the LP relaxation before solving the MILP one and then, we use the dual information given by the LP relaxation to define the $w(j, J)$ weights. Again, there is a trade-off between spending more time solving the LP relaxations and selecting better branching variables using dual information.

\begin{table}[!htpb]
\centering
{\begin{tabular}{l|lrrr}
\toprule
     Test sets &               &  \multicolumn{1}{c}{{\sc sum}} & \multicolumn{1}{c}{{\sc range}}& \multicolumn{1}{c}{{\sc dual}} \\
\midrule
					& Unsolved \nf{(26/139)} & 26\phantom{.00} & 37\phantom{.00}\phantom{\%} & 26\phantom{.00}\phantom{\%}  \\
					&$\text{Gap} = \infty$ \nf{(17)} & 19\phantom{.00} & 27\phantom{.00}\phantom{\%} & 19\phantom{.00}\phantom{\%}  \\
					&$\text{UB} = \infty$ \nf{(17)} & 19\phantom{.00} & 27\phantom{.00}\phantom{\%} & 19\phantom{.00}\phantom{\%} \\
MINLPLIB	&$\text{LB} = -\infty$ \nf{(0)} & 0\phantom{.00} & 0\phantom{.00}\phantom{\%} & 0\phantom{.00}\phantom{\%}\Bstrut \\ \cdashline{2-5}
					& Gap \nf{(6)} & 0.79 & +5.04\% & +0.10\%\Tstrut  \\
					& Time \nf{(37)} & 71.44 & +122.73\% & $-1.57\%$  \\
					& Nodes \nf{(102)} & 6.56 & +21.72\% & $-0.15\%$  \\
			\midrule
					& Unsolved \nf{(64/139)} & 65\phantom{.00} & 66\phantom{.00}\phantom{\%} & 64\phantom{.00}\phantom{\%}  \\
					&$\text{Gap} = \infty$ \nf{(25)} & 26\phantom{.00} & 27\phantom{.00}\phantom{\%} & 25\phantom{.00}\phantom{\%}  \\
					&$\text{UB} = \infty$ \nf{(25)} & 26\phantom{.00} & 27\phantom{.00}\phantom{\%} & 25\phantom{.00}\phantom{\%} \\
QPLIB			&$\text{LB} = -\infty$ \nf{(0)} & 0\phantom{.00} & 0\phantom{.00}\phantom{\%} & 0\phantom{.00}\phantom{\%}\Bstrut \\ \cdashline{2-5}
					& Gap \nf{(39)} & 9.14 & +5.97\% & +0.06\%\Tstrut \\
					& Time \nf{(75)} & 241.42 & +0.53\% & $-2.35$\%  \\
					& Nodes \nf{(73)} & 1.10 & +0.39\% & $-0.38$\%  \\
\bottomrule
\end{tabular}}
\caption{Comparison between different branching rules.}\label{table:rules}
\end{table}

Table~\ref{table:rules} shows the performance associated to the three aforementioned branching rules. The first thing we observe is that {\sc range} is the rule delivering the worst performance. This suggests that, despite being good rule for continuous problems, it may not be so informative for problems with binary variables. Interestingly, {\sc dual} slightly outperforms {\sc sum}. We see that in MINLPLIB, the worsening in the optimality gap is only $0.1$\% (over 6 instances), while the improvement on running times is $1.57$\% (over 37 instances), which is partly due to a $0.15\%$ improvement in Nodes. The results in QPLIB are qualitatively very similar but {\sc dual} managed to reduce in one unit both Unsolved and $\text{Gap}=\infty$. Thus, in the following sections we consider this new approach as the default branching rule.

\subsection{Improving upper bound through nonlinear local solvers}\label{sec:ubs}
Despite our efforts in the previous sections, our RLT-based algorithm still has the weakness of not being able to obtain an upper bound in a significant number of instances. The RLT-based algorithm, as introduced in \cite{Sherali1992}, obtains an upper bound whenever the solution of the LP relaxation at a given nodes returns a solution feasible to~\ref{eq:PO}. A limitation of this approach is that it cannot obtain an upper bound until a branch is closed, which may only occur after a large number of iterations. \cite{Dalkiran2016} study the standard approach of computing upper bounds for \ref{eq:PO} by relying on calls to an auxiliary nonlinear local solver. The results in \cite{GonzalezRodriguez2023} also confirm the significant impact of these calls on the performance of the RLT-based algorithm. In this section we explore different approaches involving not only calls to an NLP local solver, but also calls to an MINLP one, which are computationally more demanding.

Calls to local solvers involve trade-offs and, since finding the right frequencies to call them is not straightforward, global solvers are normally calibrated experimentally. Interestingly, as discussed in \cite{Khajavirad2018}, in BARON \citep{baron}, to ``increase the likelihood of finding better upper bounds on the global solution, we devised a dynamic local solver selection scheme that alternates among various solvers based on their performance in the branch-and-bound tree''. In the analysis below we present a series of much simpler (and static) strategies.

\subsubsection{Mixed-integer nonlinear local solver}
In this section we study the impact of using an MINLP local solver to increase the number of instances for which the RLT-based algorithm returns upper bounds. Since calls to a MINLP local solver are computationally demanding, we try to resort to it a reduced number of times, even just one or two calls in the whole branch-and-bound tree. In the context of solving MILP problems, there is some agreement on the fact that getting good upper bounds at the beginning of the algorithm helps to prune branches, reducing the size of the branch-and-bound tree. In this study we explore alternative strategies, including the relatively novel idea of waiting until the end to call the (computationally costly) MINLP local solver, and doing so only  if no upper bound has been previously found by the algorithm. The advantage of such an approach is that we only spend the extra computational time when it is really needed.

We argue now why, although the above approach might delay pruning, this should not negatively impact the performance of the algorithm. To illustrate, consider the baseline RLT-based algorithm in Figure~\ref{fig:RLTalg} and let $v[\text{\ref{eq:PO}}]$ be the optimal value of~\ref{eq:PO}. Then, the key observation is that the algorithm will never solve a subproblem~$k$ whose potential $LB^k$ is larger than $v[\text{\ref{eq:PO}}]$: if at some point of the branch-and-bound algorithm UB is strictly larger than $v[\text{\ref{eq:PO}}]$ (even $+\infty$), then there has to be at least one subproblem~$l$ in the queue with potential $LB^l\leq v[\text{\ref{eq:PO}}]$, so no subproblem~$k$ with $LB^k>v[\text{\ref{eq:PO}}]$ will be chosen to be solved. Of course, having good upper bounds early on will allow to prune these problems, but they will not be solved anyway. Suppose, in particular, that we follow the new approach of calling the MINLP solver at the end, and that the algorithm reaches the final phase without having obtained a valid upper bound. Looking at the scheme in Figure~\ref{fig:RLTalg}, this implies that no pruning has taken place made during the course of the algorithm and that all pruning will be done ``at once'', provided that the call to the MINLP local solver delivers a good upper bound.

Despite the above arguments, the approach of waiting until the end to call the MINLP local solver might have limitations in some cases. First, the size of the queue~$Q$ might become very large, so one has to anticipate potential memory issues. Second, there might be instances for which, at some point of the branch-and-bound tree, $LB$ is equal to $v[\text{\ref{eq:PO}}]$ or very close to it. If moreover, there are many subproblems in the queue whose lower bound is equal to  $v[\text{\ref{eq:PO}}]$ or very close to it, getting $UB=v[\text{\ref{eq:PO}}]$ would immediately terminate the algorithm but, on the other hand, not having an upper bound might lead to unnecessarily solving those problems with potential close to $v[\text{\ref{eq:PO}}]$. In order to prevent such a situation from happening, we have implemented another approach in which the MINLP local solver is called not only at the end of the algorithm, but also if it identifies that LB has been stuck for a while. Finally, recall that, as discussed in Section~\ref{sec:raposa}, the selection of the branching point in \solver{RAPOSa}~4.2.0 is such that it branches on the value of the best available solution, $\x^{best}_p$, whenever possible, so having better upper bounds early on may help to improve branching point selection.

In Table~\ref{table:knitro} we present the following approaches for calling the MINLP local solver. The first column is the best configuration from the previous section, with no calls to the MINLP local solver. In the second column, we call the MINLP local solver at the beginning of the algorithm a with a time limit of 5\% of the total running time, \emph{i.e.}, 180 seconds out of 3600 seconds. In the third column we do not call a MINLP local solver at the beginning, but we reserve 180 seconds at the end of the execution. If we reach 3420 seconds without getting an upper bound, we then call the MINLP local solver. Finally, in the last column, we adopt the same approach as the third column but introduce an additional call to the MINLP local solver if the lower bound becomes stuck,
regardless of whether we already have an upper bound. We consider that the lower bound becomes stuck if in the last $\sqrt{2N}+1$ iterations the relative change on the lower bound is below $0.001$, where $N$ is the number of variables of the underlying problem.\footnote{We considered different variations of the definition of ``becoming stuck'', with very similar results.} Furthermore, we add to Table~\ref{table:knitro} one new row (MINLP-LS time) to represent the mean time used by the MINLP local solver. We use the standard mean instead of the geometric one because there are instances in which the time used by the MINLP local solver is $0$. We also put absolute values instead of percentages because the value in the reference version is 0.

\begin{table}[!htpb]
\centering
{\begin{tabular}{l|lrrrr}
\toprule
	\multicolumn{1}{c}{}					 &               &  \multicolumn{1}{c}{MINLP-LS} & \multicolumn{1}{c}{MINLP-LS}& \multicolumn{1}{c}{MINLP-LS}& \multicolumn{1}{c}{MINLP-LS} \\
	Test sets 					 &               &  \multicolumn{1}{c}{never} & \multicolumn{1}{c}{beginning}& \multicolumn{1}{c}{end}& \multicolumn{1}{c}{end+stuck} \\ 
\midrule
						& Unsolved \nf{(25/139)} & 26\phantom{.00} & 28\phantom{.00}\phantom{\%} & 27\phantom{.00}\phantom{\%} & 27\phantom{.00}\phantom{\%} \\
						&$\text{Gap} = \infty$ \nf{(2)} & 19\phantom{.00} & 3\phantom{.00}\phantom{\%} & 3\phantom{.00}\phantom{\%} & 3\phantom{.00}\phantom{\%} \\
						&$\text{UB} = \infty$ \nf{(2)} & 19\phantom{.00} & 3\phantom{.00}\phantom{\%} & 3\phantom{.00}\phantom{\%} & 3\phantom{.00}\phantom{\%}\\
MINLPLIB		&$\text{LB} = -\infty$ \nf{(0)} & 0\phantom{.00} & 0\phantom{.00}\phantom{\%} & 0\phantom{.00}\phantom{\%} & 0\phantom{.00}\phantom{\%}\Bstrut \\ \cdashline{2-6}
						& Gap \nf{(9)} & 0.21 & +215.03\%&+204.80\%	&+188.20\%\Tstrut \\
						& Time \nf{(49)} & 31.32 & +107.27\%&+19.43\%&	+20.86\%\\
						& Nodes \nf{(111)} & 9.29 & $-$17.39\%&$-$3.19\%&$-$3.16\% \\
						& MINLP-LS time \nf{(139)} & 0.00 & 20.80\phantom{\%} & 15.58\phantom{\%} & 18.25\phantom{\%}\\						
				\midrule
						& Unsolved \nf{(63/139)} & 64\phantom{.00} & 63\phantom{.00}\phantom{\%} & 64\phantom{.00}\phantom{\%} & 64\phantom{.00}\phantom{\%}\\
					  & $\text{Gap} = \infty$ \nf{(5)} & 25\phantom{.00} & 6\phantom{.00}\phantom{\%} & 7\phantom{.00}\phantom{\%} & 7\phantom{.00}\phantom{\%} \\
						&$\text{UB} = \infty$ \nf{(5)} & 25\phantom{.00} & 6\phantom{.00}\phantom{\%} & 7\phantom{.00}\phantom{\%} & 7\phantom{.00}\phantom{\%}\\
QPLIB				&$\text{LB} = -\infty$ \nf{(0)} & 0\phantom{.00} & 0\phantom{.00}\phantom{\%} & 0\phantom{.00}\phantom{\%} & 0\phantom{.00}\phantom{\%}\Bstrut \\ \cdashline{2-6}
						& Gap \nf{(40)} & 7.30 & +0.32\%&$-$0.65\%&$-$1.08\%\Tstrut\\
						& Time \nf{(76)} & 243.54 & +16.28\%&$-$2.09\%&$-$1.78\% \\
						& Nodes \nf{(74)} & 1.36 & $-$0.30\%&$-0.66$\%&$-$0.66\%\\
						& MINLP-LS time \nf{(139)} & 0.00 & 35.42\phantom{\%} & 22.86\phantom{\%} & 24.74\phantom{\%}\\						
\bottomrule
\end{tabular}}
\caption{Comparison between different approaches of calling the MINLP local solver.}\label{table:knitro}
\end{table}

Table~\ref{table:knitro} shows that calling a MINLP local solver really helps the RLT-based algorithm to find more upper bounds. Specifically, calling a MINLP local solver reduces the number of instances in which the RLT-based algorithm fails to return an upper bound from 19 to 3 in MINLPLIB and from 35 to 6-7 in QPLIB. Moreover, the last two approaches exhibit superior performance in terms of optimality gaps and running times compared to the second one. This improvement is due to the reduction in time consumption by the MINLP local solver, as it is called fewer times. Thus, calling the MINLP solver at the end of the algorithm seems to outperform the more conventional approach of calling it at the beginning to obtain upper bounds more quickly. Note that, in MINLPLIB, calling the MINLP local solver at the beginning seems to notably reduce Nodes in the solved instances, which may be the result of improved branching point selection, given the availability of better upper bounds early on, as discussed above. However, the additional time spent by the MINLP solver offsets this advantage. Finally, considering the last two approaches, they exhibit a very similar behavior. Given that the last one might provide some additional robustness, we take it as the reference version for the following sections. 

\subsubsection{Continuous nonlinear local solver}

Given that the calls to the MINLP local solver are computationally expensive, in this section we try to obtain more upper bounds with the NLP local solver, which is computationally cheaper. If we are successful at getting more upper bounds with the NLP local solver, then we can spare the call to the MINLP local solver at the end of the algorithm. Following \cite{GonzalezRodriguez2023}, the default version of \solver{RAPOSa} for continuous problems calls the NLP local solver at the root node and whenever the iteration count in the branch-and-bound algorithm is a power of~2. At the root node no initial solution is provided to the NLP solver and in the rest of nodes the initial solution is just the optimal solution of the LP relaxation at the corresponding node. They mention that they tried other strategies, including different frequencies for the NLP calls, but they did not observe a significant impact on performance.

As briefly discussed in Section~\ref{sec:raposa}, when dealing with mixed-integer polynomial optimization problems we have so far relied on a straightforward adaptation of the above approach, in which we take the solution provided by the NLP local solver and round the integer variables. If this rounded solution is feasible to~\ref{eq:PO}, then we have an upper bound. However, since in this mixed-integer setting this approach is likely to be less effective at returning upper bounds than the approach for the continuous one, we have increased the frequency at which the NLP local solver is called. Specifically, we call it at the root node and whenever the total number of solved nodes in the branch-and-bound tree is a power of~$1.5$. This is the approach we have used so far in all the numeric analyses.

In this section we study the performance of three different approaches to call the NLP local solver. The first one, NLP-LS-round, is the one used so far, in which the initial solution for the NLP local solver is the optimal solution at the current node and then the integer variables in the local optimum returned by the NLP local solver are rounded. The second approach, NLP-LS-fix, consists of fixing the integer variables at their values in the optimal solution of the current node before passing it to the NLP local solver. The main advantage of this approach is that, if the NLP local solver returns a solution, it will be feasible to~\ref{eq:PO} but, at the same time, in most calls to the NLP solver it will be given infeasible subproblems. The third approach, NLP-LS-round+fix, sequentially performs the two previous approaches. Regarding the quality of the upper bounds, by definition, NLP-LS-round+fix always returns something at least as good as NLP-LS-round and NLP-LS-fix. As for computational times, NLP-LS-round+fix always consumes more time than NLP-LS-round and NLP-LS-fix, while NLP-LS-fix should be in general faster than NLP-LS-round, since some variables are fixed.\footnote{Note that NLP-LS-fix needs the optimal solution of the current node so it cannot be used before solving the root node. Thus, we always do the first call to the MINLP local solver using NLP-LS-round and then, after solving the root node, we start using the corresponding approach.}

\begin{table}[!htpb]
\centering
{\begin{tabular}{l|lrrr}
\toprule
	\multicolumn{1}{c}{}	 &       &  \multicolumn{1}{c}{NLP-LS} & \multicolumn{1}{c}{NLP-LS}& \multicolumn{1}{c}{NLP-LS}\\
	Test sets 					 &               &  \multicolumn{1}{c}{round} & \multicolumn{1}{c}{fix}& \multicolumn{1}{c}{round+fix}\\   
\midrule
						& Unsolved \nf{(26/139)} & 27\phantom{.00} & 26\phantom{.00}\phantom{\%} & 26\phantom{.00}\phantom{\%} \\
						&$\text{Gap} = \infty$ \nf{(2)} & 3\phantom{.00} & 2\phantom{.00}\phantom{\%} & 2\phantom{.00}\phantom{\%} \\
						&$\text{UB} = \infty$ \nf{(2)} & 3\phantom{.00} & 2\phantom{.00}\phantom{\%} & 2\phantom{.00}\phantom{\%} \\
MINLPLIB		&$\text{LB} = -\infty$ \nf{(0)} & 0\phantom{.00} & 0\phantom{.00}\phantom{\%} & 0\phantom{.00}\phantom{\%}\Bstrut \\ \cdashline{2-5}
						& Gap \nf{(24)} & 0.90 & $-16.44$\% & $-16.46$\%\Tstrut \\
						& Time \nf{(37)} & 60.76 & $-11.83$\% & $-8.70$\% \\
						& Nodes \nf{(112)} & 9.12 &	$-6.23$\%	&$-6.23$\%\\ 
						&MINLP-LS time \nf{(139)} & 18.25 & 15.29\phantom{\%} & 15.28\phantom{\%} \\
						&NLP-LS time \nf{(139)} & 1.75 &\phantom{\%} 0.92\phantom{\%} & 2.00\phantom{\%} \\						
				\midrule
						& Unsolved \nf{(63/139)} & 64\phantom{.00} & 63\phantom{.00}\phantom{\%} & 63\phantom{.00}\phantom{\%} \\
						&$\text{Gap} = \infty$ \nf{(5)} & 7\phantom{.00} & 6\phantom{.00}\phantom{\%} & 5\phantom{.00}\phantom{\%} \\
						&$\text{UB} = \infty$ \nf{(5)} & 7\phantom{.00} & 6\phantom{.00}\phantom{\%} & 5\phantom{.00}\phantom{\%} \\
QPLIB				&$\text{LB} = -\infty$ \nf{(0)} & 0\phantom{.00} & 0\phantom{.00}\phantom{\%} & 0\phantom{.00}\phantom{\%}\Bstrut \\ \cdashline{2-5}
						& Gap \nf{(57)} &3.26 & +6.00\%	&$+$6.04\%\Tstrut \\
						& Time \nf{(76)} & 239.19 & $-$4.09\%	&$-$3.47\% \\
						& Nodes \nf{(75)} & 1.47 & $+$0.61\%	&$+$0.61\% \\
						&MINLP-LS time \nf{(139)} & 24.74 & 25.34\phantom{\%} & 25.37\phantom{\%} \\
						&NLP-LS time \nf{(139)} & 12.16 & 8.47\phantom{\%} & 13.13\phantom{\%} \\						
\bottomrule
\end{tabular}}
\caption{Comparison between different approaches for the calls to the NLP local solver.}\label{table:ipopt_approach}
\end{table}

In Table~\ref{table:ipopt_approach} we present a summary of results. We include a new row called NLP-LS time to represent the average of the time used by the NLP local solver, and keep the MINLP-LS row from Table~\ref{table:knitro} (again standard averages are used to report these times). First, we see that NLP-LS-fix improves the overall performance with respect to NLP-LS-round. Indeed, NLP-LS-fix reduces Time in both test sets and, regarding gap, it goes down by $16.44\%$ in MINLPLIB and goes up by just $6\%$ in QPLIB. Moreover, NLP-LS-fix solves one more instance in each test set and returns one additional upper bound in each test set. If we look at NLP-LS time, we see that it is smaller for NLP-LS-fix, as expected, since we are fixing the integer variables. Interestingly, NLP-LS-fix also results in a smaller value for MINLP time in MINLPLIB, probably because this approach allows to dispense with some calls to the MINLP local solver at the end of the algorithm. When comparing NLP-LS-fix and NLP-LS-round+fix, we see that they exhibit a very similar performance. NLP-LS-round+fix is slightly worse in terms of Gap and Time, probably due to the extra NLP-LS time. Yet, NLP-LS-round+fix returns one more upper bound in QPLIB. Since our ultimate goal is to return more upper bounds without compromising much the performance of the RLT-based algorithm, we choose NLP-LS-round+fix as the reference configuration.

It is worth noting that we also tested a significant number of alternative approaches regarding the NLP local solver, but they were not as competitive as the ones in Table~\ref{table:ipopt_approach}. For the sake of completeness, we briefly mention two of them. In the first one, whenever the NLP solver was called, it was done several times, with different initial solutions. These  solutions were obtained by asking the MILP solver to return a pool of feasible solutions instead of just the optimal one. The second one was simply to play around with the frequency of the calls to the NLP solver, being more aggressive by calling it at powers of $p$, with different values for $p$, smaller than $1.5$.

\subsection{Comparison with other state-of-the-art solvers}

The different approaches presented in the previous sections have led to an RLT-based algorithm that is much more competitive at finding both lower and upper bounds. Once this goal has been achieved, we revisit again the comparison with the state-of-the-art solver \solver{BARON}, using a version of \solver{RAPOSa} equipped with these enhancements. This time we omit \solver{Couenne} from the comparison, given its subpar performance in the results in Section~\ref{sec:othersolvers}.
On the other hand, we keep the version of \solver{RAPOSa} used in Section~\ref{sec:othersolvers}, in order to have a reference on which to assess the impact of the enhancements. The main goal in the previous sections was to be able to improve in metrics $\text{Gap}=\infty$, $\text{UB}=\infty$, and $\text{LB}=\infty$ without compromising too much the performance of the RLT-based algorithm in Gap and Time. Again, the executions have been run with a common limit of 8 threads to all solvers, which results in slightly worse results for \solver{RAPOSa} than the ones reported in the previous sections, were it had been run letting it use the 32 cores available. Also, given that \solver{Couenne} is not present, the problems excluded for the computation of each metric are slightly different from those excluded in Table~\ref{table:solvers1}.

\begin{table}[!htpb]
\centering
{\begin{tabular}{l|lrrrr}
\toprule
     Test sets &          &  \solver{RAPOSa} old &  \solver{RAPOSa} new & \solver{BARON} (\solver{CBC})& \solver{BARON} (\solver{CPLEX})  \\
\midrule
					& Unsolved \nf{(16/139)} & 30\phantom{.000}& 29\phantom{.00}\phantom{\%} &30\phantom{.00}\phantom{\%}&	25\phantom{.00}\phantom{\%} \\
					&$\text{Gap} = \infty$ \nf{(0)} & 22\phantom{.000} & 3\phantom{.00}\phantom{\%} & 0\phantom{.00}\phantom{\%} & 1\phantom{.00}\phantom{\%} \\
					&$\text{UB} = \infty$ \nf{(0)} & 22\phantom{.000} & 3\phantom{.00}\phantom{\%}& 0\phantom{.00}\phantom{\%} & 1\phantom{.00}\phantom{\%}  \\
MINLPLIB	&$\text{LB} = -\infty$ \nf{(0)} & 7\phantom{.000} & 0\phantom{.00}\phantom{\%}& 0\phantom{.00}\phantom{\%} & 0\phantom{.00}\phantom{\%}\Bstrut \\ \cdashline{2-6}
					& Gap \nf{(20)} & 0.011& +56.66\%& +517.48\% & +179.97\%\Tstrut  \\
					& Time \nf{(52)} & 87.597& +23.42\% & +23.18\% & $-81.31\%$\\
					\midrule
					& Unsolved \nf{(55/139)} & 67\phantom{.000} & 68\phantom{.00}\phantom{\%}& 105\phantom{.00}\phantom{\%}	& 79\phantom{.00}\phantom{\%}	 \\
					&$\text{Gap} = \infty$ \nf{(0)} & 66\phantom{.000}&7\phantom{.00}\phantom{\%} &0\phantom{.00}\phantom{\%}&	0\phantom{.00}\phantom{\%}	 \\
					&$\text{UB} = \infty$ \nf{(0)} & 27\phantom{.000}&	7\phantom{.00}\phantom{\%} &0\phantom{.00}\phantom{\%}	&0\phantom{.00}\phantom{\%} \\
QPLIB			&$\text{LB} = -\infty$ \nf{(0)} & 55\phantom{.000} & 0\phantom{.00}\phantom{\%}& 0\phantom{.00}\phantom{\%} & 0\phantom{.00}\phantom{\%}\Bstrut \\ \cdashline{2-6}
					& Gap \nf{(53)} & 0.001 & +17.19\% & +10396.22\% & +1892.00\%\Tstrut\\
					& Time \nf{(84)} & 396.682& $-1.80\%$ & +241.57\% & $-22.72\%$  \\
\bottomrule
\end{tabular}}
\caption{Comparison of two configurations of \solver{RAPOSa} with \solver{BARON}.}\label{table:solvers2}
\end{table}

The results in Table~\ref{table:solvers2} show that the new version of \solver{RAPOSa} clearly accomplishes the desired goals. For instance, if we look at MINLPLIB, we have that $\text{Gap}=\infty$ goes down from~22 to~3, and it comes at the relatively small price of increasing Gap from $0.011$ to $0.0176$ and Time from $87.597$ to $108.111$. Also, note that the increase in Gap is computed on just 20 instances, roughly the same number as the new instances in which \solver{RAPOSa} is now capable of returning a finite gap. The results for QPLIB are qualitatively similar, but even better. $\text{Gap}=\infty$ goes down from~66 to~7, with Gap increasing just $17\%$, from $0.0012$ to $0.0014$, and Time slightly improving from $396.682$ to	$389.525$. The increase in Gap is computed on 53 instances, a number again comparable to the new instances in which \solver{RAPOSa} is now capable of returning a finite gap, 59. If we move now to the comparison with \solver{BARON}, we see that the new version of \solver{RAPOSa} is now quite competitive in all metrics. \solver{RAPOSa} is still slightly behind at finding upper bounds, but it is now much closer to \solver{BARON} (\solver{CPLEX}) than it was before. Overall, \solver{RAPOSa} seems to be superior in QPLIB and \solver{BARON} (\solver{CPLEX}) at MINLPLIB and, again, in both test sets \solver{RAPOSa} is better at closing the gap in difficult instances and \solver{BARON} (\solver{CPLEX}) is better at solving the not-so-difficult ones.

\begin{figure}[!htpb]
\centering
\begin{subfigure}{0.48\textwidth}
\includesvg[pretex=\scriptsize,width=\textwidth]{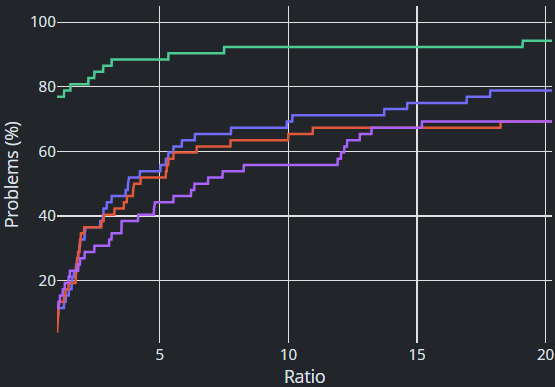}
    \caption{Running time in MINLPLIB.}
    \label{fig:timeminlpinit}
\end{subfigure}
\hfill
\begin{subfigure}{0.48\textwidth}
    \begin{tikzpicture}
     \node (fig1) at (0,0){\includesvg[pretex=\scriptsize,width=\textwidth]{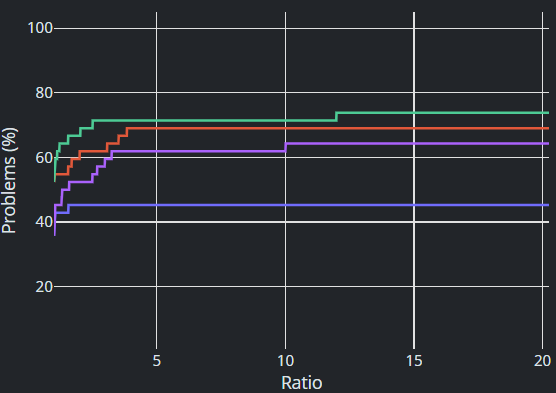}};
     \node (fig2) at (3.1,-0.9){\includesvg[pretex=\scriptsize,width=0.5\textwidth]{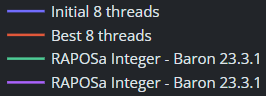}}; 
    \end{tikzpicture}
    \caption{Optimality gap in MINLPLIB.}
    \label{fig:gapminlpfinal}
\end{subfigure}

\begin{subfigure}{0.48\textwidth}
    \includesvg[pretex=\scriptsize,width=\textwidth]{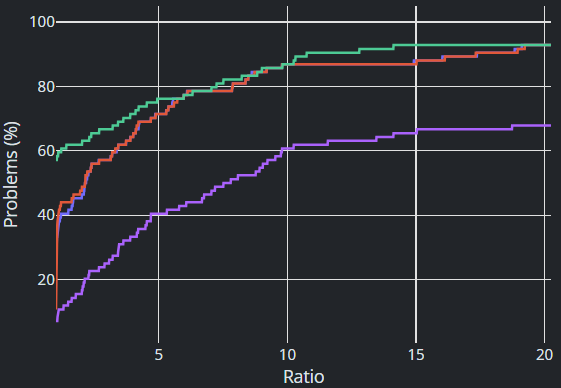}
    \caption{Running time in QPLIB.}
    \label{fig:timeqplibfinal}
\end{subfigure}
\hfill
\begin{subfigure}{0.48\textwidth}
    \includesvg[pretex=\scriptsize,width=\textwidth]{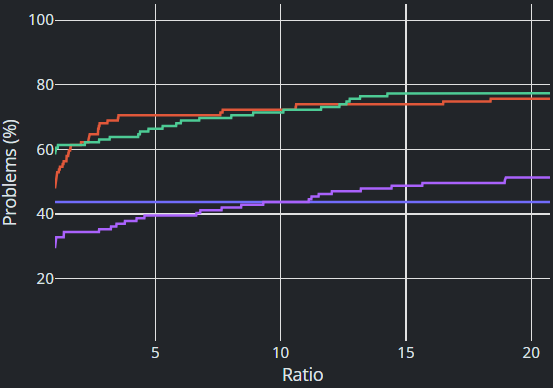}
    \caption{Optimality gap in QPLIB.}
    \label{fig:gapqplibfinal}
\end{subfigure}
\caption{Performance profiles comparing two configurations of \solver{RAPOSa} with \solver{BARON}.}
\label{fig:solvers2}
\end{figure}

To conclude the numerical analysis, in Figure~\ref{fig:solvers2} we present a last set of performance profiles, which are again consistent with those in Table~\ref{table:solvers2}, and allow to see quite clearly the improvements obtained with the enhancements in \solver{RAPOSa} and the associated trade-offs. Regarding running times, \solver{BARON} (\solver{CPLEX}) is ahead in both test sets, and by a big margin in MINLPLIB. Also, the two versions of \solver{RAPOSa} perform very similarly in QPLIB, while the new version is slightly worse in MINLPLIB. The big step forward in the new version can be appreciated when looking at the optimality gaps, where the new version of RAPOSa is now very competitive with \solver{BARON} (\solver{CPLEX}), and particularly so in QPLIB.

\section{Conclusions}\label{sec:conclusions}

In this paper we have extended a continuous RLT-based solver for polynomial optimization problems to handle mixed-integer problems. We have identified the main weaknesses of a first implementation of this new algorithm, and then introduced some enhancements that have successfully mitigated them, leading to a version of \solver{RAPOSa} that is highly competitive with the state-of-the-art solver \solver{BARON}. Interestingly, one of these enhancements consists of calling the auxiliary MINLP local solver at the end of the algorithm, which, to the best of our knowledge, is not the standard practice, and has shown to deliver very good results in our setting. It might be worth exploring the impact of similar approaches for other well established global optimization solvers.

Given the promising results obtained in this paper, it is clear that RLT-based algorithms can be very effective at solving (nonconvex) MINLP polynomial optimization problems. Yet, there is significant room for improvement in the implementation described in this paper, and more research is needed to fully understand the potential of RLT at solving these general problems. In particular, this first implementation heavily relies on the tightness provided by the bound-factor constraints, whose caveat is that their number grows exponentially on the size of the problem. These bound-factor constraints are defined for all the variables of the problem, regardless of whether they are continuous or not. A natural direction of future research is to rely on alternative reformulations for the binary variables along the lines of \cite{elloumi2023}, which can reduce the degree of the reformulated problem and, therefore, the number of bound-factor constraints. Since such reformulations would lead to looser relaxations, one might then build upon the literature on valid inequalities to recover some of the lost tightness, as in \cite{delpiaSahinidis2020}.

\section*{Acknowledgments}
Project PID2020-116587GB-I00 supported by MINECO through the ERDF and project PID2021-124030NB-C32 funded by MICIU/AEI/10.13039/501100011033/ and by ERDF/EU. This research was also funded by Grupos de Referencia Competitiva ED431C-2021/24 from the Consellería de Cultura, Educación e Universidades, Xunta de Galicia. Iria Rodríguez-Acevedo acknowledges the support from Consellería de
Cultura, Educación, Formación Profesional e Universidades, Xunta de Galicia, through grant ED481A-2023-061.

\bibliographystyle{ecta}
\bibliography{references}





\end{document}